\input amstex.tex
\input epsf
\documentstyle{amsppt}
\magnification=1200
\baselineskip=13pt
\hsize=6.5truein
\vsize=8.9truein
\parindent=20pt
\newcount\sectionnumber
\newcount\equationnumber
\newcount\thnumber
\newcount\countrefno

\def\ifundefined#1{\expandafter\ifx\csname#1\endcsname\relax}
\def\assignnumber#1#2{%
	\ifundefined{#1}\relax\else\message{#1 already defined}\fi
	\expandafter\xdef\csname#1\endcsname
        {\the\sectionnumber.\the#2}}%
%
%
\def\newsec{
  \global\advance\sectionnumber by 1
  \global\equationnumber=0
  \global\thnumber=0
  \the\sectionnumber .\ }
%
\def\eq#1{\relax
  \global\advance\equationnumber by 1
  \assignnumber{EN#1}\equationnumber
  {\rm \csname EN#1\endcsname}}
\def\eqtag#1{\ifundefined{EN#1}\message{EN#1 undefined}{\sl (#1)}%
  \else\thetag{\csname EN#1\endcsname}\fi}
%
%
\def\thname#1{\relax
  \global\advance\thnumber by 1
  \assignnumber{TH#1}\thnumber
  \csname TH#1\endcsname}
\def\thtag#1{\ifundefined{TH#1}\message{TH#1 undefined}{\sl #1}%
  \else\csname TH#1\endcsname\fi}
%
\comment
\def\eq{}
\def\eqtag#1{(#1)}
\def\thname{}
\def\thtag{}
\endcomment
%
%
\global\countrefno=1

\def\refno#1{\xdef#1{{\the\countrefno}}
\global\advance\countrefno by 1}
\def\R{{\Bbb R}}
\def\N{{\Bbb N}}
\def\C{{\Bbb C}}

\def\Z{{\Bbb Z}}
\def\T{{\Bbb T}}

\def\Zp{{{\Bbb Z}_{\geq 0}}}

\def\hf{{1\over 2}}

\def\al{\alpha}
\def\be{\beta}
\def\ga{\gamma}
\def\Ga{\Gamma}
\def\de{\delta}
\def\De{\Delta}

\def\ep{\varepsilon}

\def\th{\theta}
\def\la{\lambda}
\def\vp{\varphi}

\refno{\AndrA}
\refno{\AskeW}
\refno{\BoneCGST}
\refno{\BrowEI}
\refno{\BustS}
\refno{\CiccKK}
\refno{\GaspR}
\refno{\GuptIM}
\refno{\Isma}
\refno{\IsmaMS}
\refno{\IsmaR}
\refno{\Kake}
\refno{\KakeMU}
\refno{\KoekS}
\refno{\KoelPhD}
\refno{\KoelDMJ}
\refno{\KoelIM}
\refno{\KoelJCAM}
\refno{\KoelAAM}
\refno{\KoelJAT}
\refno{\KoelSbig}
\refno{\KoelSsu}
\refno{\KoelSAW}
\refno{\KoorJF}
\refno{\KoorLNM}
\refno{\KoorSIAM}
\refno{\KoorLNSFQG}
\refno{\KoorS}
\refno{\MasuMNNSU}
\refno{\Moak}
\refno{\Noum}
\refno{\NoumM}
\refno{\Rose}
\refno{\StokK}
\refno{\Susl}
\refno{\SuslPP}
\refno{\VaksKmotion}
\refno{\VaksKsu}
\refno{\Wats}
\topmatter
\title The Askey-Wilson function transform scheme\endtitle
\author Erik Koelink and Jasper V.~Stokman 
\endauthor
\address Technische Universiteit Delft, Faculteit
Informatietechnologie en Systemen, Afd. Toegepaste Wiskundige
Analyse, Postbus 5031, 2600 GA Delft, the Netherlands\endaddress
\email koelink\@twi.tudelft.nl\endemail
\address Centre de Math\'ematiques de Jussieu, Universit\'e 
Paris 6 Pierre et Marie Curie, 4, place Jussieu, 
F-75252 Paris Cedex 05, France\endaddress
\email stokman\@math.jussieu.fr\endemail
\date December 23, 1999\enddate 
\abstract  In this paper we present an addition to Askey's
scheme of $q$-hypergeometric orthogonal polynomials 
involving classes of $q$-special functions which do not consist
of polynomials only. The special functions are $q$-analogues of the
Jacobi and Bessel function. The generalised orthogonality
relations and the second order $q$-difference
equations for these families are given. Limit transitions
between these families are discussed. The quantum group theoretic
interpretations are discussed shortly. 
\endabstract
\endtopmatter

\document

\head \newsec Introduction
\endhead

The Askey-scheme of hypergeometric orthogonal polynomials is
a scheme containing various known sets of orthogonal polynomials
that can be written in terms of hypergeometric series, see 
e.g. Askey and Wilson \cite{\AskeW}, Koekoek and Swarttouw 
\cite{\KoekS}, 
Koornwinder \cite{\KoorLNSFQG}.
A typical entry is the set of Jacobi polynomials defined by
$$
R_n^{(\al,\be)}(x)= {}_2F_1\left( {{-n,n+\al+\be+1}\atop
{\al+1}}; {{1-x}\over{2}}\right),
\qquad n\in\Zp,
\tag\eq{010}
$$
where we use the standard notation for hypergeometric series,
see e.g. \cite{\GaspR}. The Jacobi polynomials 
are orthogonal with respect to the beta distribution 
$(1-x)^\al(1+x)^\be$ on the interval $[-1,1]$. 
Moreover, we have the much bigger $q$-analogue 
of the Askey-scheme having the Askey-Wilson polynomials and Racah
polynomials at the top level with four degrees of freedom (apart
from $q$), see e.g. \cite{\KoekS}.

The relation with schemes of non-polynomial special
functions are less well-advertised, and here we discuss the 
$q$-analogue of the scheme displayed in Figure~1.1. 
The Jacobi function transform is an integral transform 
on $[0,\infty)$ in which
the kernel is given by a Jacobi function 
$$
\phi^{(\al,\be)}_\la(t)= {}_2F_1\left( {{\hf(\al+\be+1-i\la),
\hf(\al+\be+1+i\la)}\atop{ \al+1}};-\sinh^2t\right) 
\tag\eq{020}
$$ 
for $|\sinh t|<1$, which has a one-valued analytic continuation
to $-\sinh^2 t\in\C\backslash [1,\infty)$, 
see Koornwinder's survey paper \cite{\KoorJF}
for a nice introduction and references. We see that we
may view the Jacobi function as an analytic continuation in 
its degree of the Jacobi polynomial; replace $n$ in \eqtag{010}
by $\hf(i\la-\al-\be-1)$ and $x$ by $\cosh 2t$ to get 
the Jacobi function of \eqtag{020}.

\topinsert
\hskip2truecm
{\epsfxsize 3.0truein \epsfysize 1.75truein \epsfbox{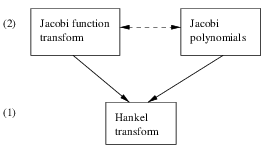}}
\botcaption{Figure 1.1}
Jacobi function scheme.
\endcaption
\endinsert 

The Jacobi function transform is then given by
$$
\aligned
g(\la) &= \int_0^\infty f(x) \phi^{(\al,\be)}_\la(t)
(2\sinh t)^{2\al+1}(2\cosh t)^{2\be+1}\, dt, \\
f(t) &= \int_0^\infty g(\la) \phi^{(\al,\be)}_\la(t) 
{{|\Ga(\hf(i\la+\al+\be+1))\Ga(\hf(i\la+\al-\be+1))|^2}\over{
4^{\al+\be+1}|\Ga(\al+1) \Ga(i\la)|^2}}\, d\la
\endaligned
\tag\eq{030}
$$
for some suitable class of functions. Here we assume that
$\al,\be\in\R$ satisfy $|\be|< \al+1$, otherwise discrete
mass points have to be added to the Plancherel measure, see 
\cite{\KoorJF, \S 2}. 

The Hankel transform is the integral transform on $[0,\infty)$ 
that has the Bessel function
$$
J_\al(x)= {{(x/2)^\al}\over{\Ga(\al+1)}}
\, {}_0F_1\left( {{-}\atop{\al+1}};-{{x^2}\over 4}\right)
\tag\eq{040}
$$
as its kernel. For suitable functions and $\al>-1$ 
the Hankel transform pair is given by,
see Watson \cite{\Wats, \S 14.3},  
$$
g(\la) = \int_0^\infty f(x)\, J_\al(x\la) x\, dx, \qquad
f(x) = \int_0^\infty g(\la)\, J_\al(x\la) \la\, d\la.
\tag\eq{050}
$$
The Hankel transform can formally be obtained as a limit case
from the orthogonality relations for the
Jacobi polynomials using the limit
$$
\lim_{N\to\infty} R_{n_N}^{(\al,\be)}(1-{{x^2}\over{2N^2}}) 
= {}_0F_1(-;\al+1;- (x\la)^2/4),\quad 
\text{$n_N/N\to \la$ as $N\to\infty$.}
\tag\eq{060}
$$
The Hankel transform
can also be viewed as a limit case of the Jacobi function 
transform by use of the limit transition
$$
\lim_{\ep\downarrow 0} \phi_{\la/\ep}^{(\al,\be)}(t\ep) = 
2^\al \Ga(\al+1) (\la t)^{-\al}\, J_\al(\la t),
\tag\eq{070}
$$
see \cite{\KoorJF, (2.34)}, and the limit transition \eqtag{070}
can be considered as a generalisation of \eqtag{060}.

\topinsert
\hskip0.2truecm
{\epsfxsize 5.75truein \epsfysize 5truein 
\epsfbox{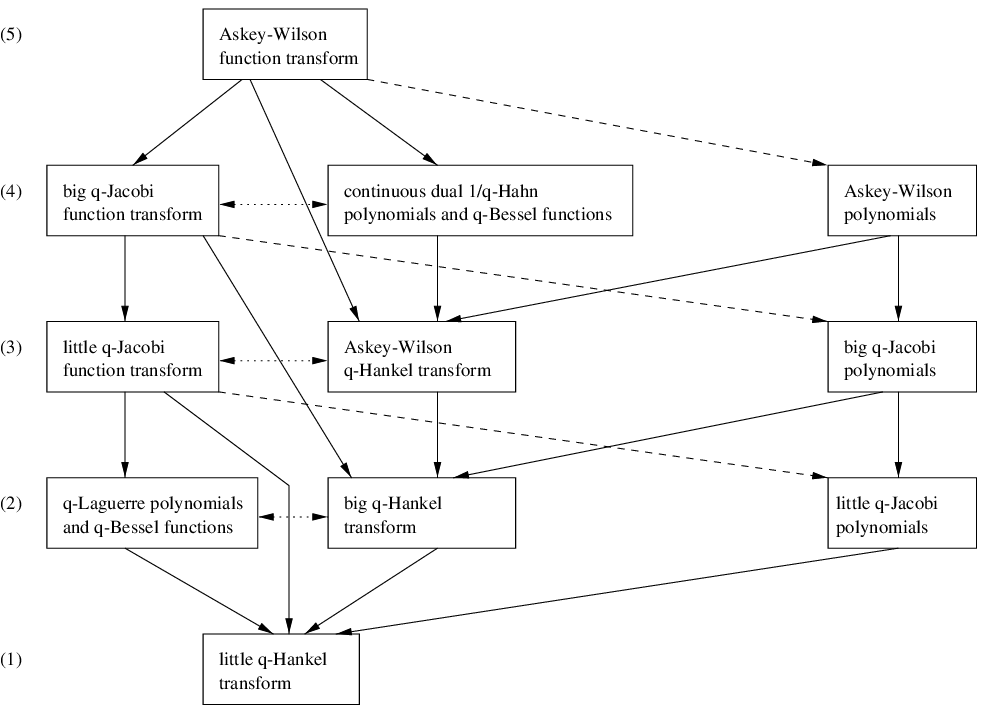}}
\botcaption{Figure 1.2}
Askey-Wilson function scheme.
\endcaption
\endinsert

So Figure~1.1 gives an addition to the Askey-scheme of
hypergeometric orthogonal polynomials by an analytic 
continuation in the spectral parameter 
(the dashed line) to the Jacobi functions, and
by a limit transition to the Bessel functions. The Jacobi functions
and polynomials have two degrees of freedom, namely $\al$ and
$\be$, and the Hankel transform has one degree of freedom,
namely the order $\al$ of the Bessel function. We consider the
Jacobi functions as the master functions, 
and the Jacobi polynomials
and the Bessel functions as derivable functions. 
We discuss the group theoretic interpretation of Figure~1.1 in
\S 7.1. 

The purpose of this paper is to give three $q$-analogues of 
Figure~1.1; one related to the Askey-Wilson polynomials, 
presented in \S 2, one related to the big $q$-Jacobi polynomials,
presented in \S 3, and one related to the little
$q$-Jacobi polynomials, presented in \S 4. This is depicted in
Figure~1.2, where the three boxes on the right hand side are
part of Askey's scheme of $q$-hypergeometric orthogonal polynomials,
see \cite{\KoekS}. These three 
$q$-analogues are related by limit transitions 
as well, and this is also depicted in Figure~1.2.
This is discussed in \S 5 and \S 6, where the boxes with the
continuous dual $q^{-1}$-Hahn polynomials and the $q$-Laguerre
polynomials are discussed.  

In Figure~1.2 the lines denote limit transitions, 
the dashed lines denote analytic continuation and the dotted
lines, which do not appear in Figure~1.1, denote that the
results of the two boxes involved are related by duality. 
In particular, we think of the Askey-Wilson function 
transform and the little $q$-Hankel transform as self-dual
transforms, i.e. the inverse transform equals the
generic Fourier transform itself, possibly for dual
parameters. 
Note that there is a striking difference between Figure~1.1
and Figure~1.2. The dashed line corresponding to 
analytic continuation in Figure~1.1 is horizontal, i.e. 
the Jacobi polynomials and the Jacobi function transform
both have 2 degrees of freedom. The dashed lines in Figure~1.2
go down, i.e. the $q$-analogue of the Jacobi function transform
has one extra degree of freedom compared to the $q$-analogue
of the Jacobi polynomial. This extra parameter is not contained in
the definition of the $q$-analogues of the Jacobi 
function, but it appears in the measure for the 
$q$-analogue of the Jacobi function transform. 

In \S 7 we discuss the (quantum) group theoretic
interpretation of Figure~1.1 and Figure~1.2, which is
also the motivation for the names for these transforms. 
We end with some concluding remarks and some open
problems. 

Let us finally note that this paper does not contain
rigorous proofs. The transform pairs are obtained 
by a spectral analysis of a second order $q$-difference
operator that is symmetric on a suitable weighted $L^2$-space, 
see the references in \S\S 2, 3 and 4.
The limit transitions are only considered on 
a formal level and are emphasised via the second order
$q$-difference operators involved.  

\demo{Notation} We use the notation for basic hypergeometric
series as in the book \cite{\GaspR} by Gasper and Rahman. 
We assume throughout $0<q<1$. 
\enddemo

\demo{Acknowledgement} The second author is 
supported by a NWO-Talent stipendium of the 
Netherlands Organization for Scientific Research (NWO).
Part of the research was done while the second author
was supported by the EC TMR network ``Algebraic Lie
Representations'', grant no. ERB FMRX-CT97-0100. 
\enddemo

\head \newsec Askey-Wilson analogue of the Jacobi function 
scheme\endhead

In this section we consider the analogue of Figure~1.1 on the
level for the Askey-Wilson case. So the second order difference
operator is
$$
Lf(x) = A(x)\bigl( f(qx)-f(x)\bigr) + 
B(x)\bigl( f(q^{-1}x)-f(x)\bigr),
\tag\eq{110}
$$
where 
$$
A(x) = {{(1-ax)(1-bx)(1-cx)(1-dx)}\over{(1-x^2)(1-qx^2)}},
\qquad B(x)=A(x^{-1}). 
\tag\eq{120}
$$
The difference operator \eqtag{110}, \eqtag{120} has been
introduced by Askey and Wilson \cite{\AskeW}. 
The general set of eigenfunctions for \eqtag{110} with
$A$ and $B$ as in \eqtag{120} is given by 
Ismail and Rahman \cite{\IsmaR}. 

\subhead \the\sectionnumber.1. The Askey-Wilson 
functions\endsubhead The 
Askey-Wilson functions are defined by 
$$
\multline
\phi_\ga(x;a;b,c;d|q)=\frac{\bigl(qax\gamma/\tilde{d}, 
qa\gamma/\tilde{d}x;q\bigr)_{\infty}}
{\bigl(\tilde{a}\tilde{b}\tilde{c}\gamma,q\gamma/\tilde{d},
q\tilde{a}/\tilde{d}, qx/d,q/dx;q\bigr)_{\infty}}\\
\times {}_8W_7\left(\tilde{a}\tilde{b}\tilde{c}\gamma/q; 
ax, a/x, \tilde{a}\gamma, \tilde{b}\gamma, \tilde{c}\gamma;
q,q/\tilde{d}\gamma\right)
\endmultline
\tag\eq{130}
$$
for $|\ga|>|q/\tilde d|$, 
where $\tilde{a}=\sqrt{q^{-1}abcd}$, $\tilde{b}=
ab/\tilde{a}$, $\tilde{c}=ac/\tilde{a}$,
and $\tilde{d}=ad/\tilde{a}$. Observe that the Askey-Wilson
function in \eqtag{130} is symmetric in $b$ and $c$, and
almost symmetric in $a$ and $b$;
$$
\phi_\ga(x;a;b,c;d|q)= \frac{(\tilde c\ga,qb/d,\tilde
c/\ga;q)_\infty}{(q\ga/\tilde d, qa/d,q/\tilde d\ga;q)_\infty}
\phi_\ga(x;b;a,c;d|q), 
\tag\eq{132}
$$
by an application of \cite{\GaspR, (III.36)}, or see 
Suslov \cite{\SuslPP}. 
Then
$L\phi_\ga = \bigl(-1-\tilde{a}^2+\tilde{a}(\ga+\ga^{-1})\bigr)
\phi_\ga$. 
Note that the invariance $x\leftrightarrow x^{-1}$ is obvious in 
\eqtag{130}. There exists a meromorphic continuation
of the Askey-Wilson function in $\ga$ which is invariant 
under $\ga\leftrightarrow \ga^{-1}$ by 
\cite{\GaspR, (III.23)}. By \cite{\GaspR, (III.23)} again we
see that the Askey-Wilson function is 
self-dual in the sense that 
$$
\phi_{\ga}(x;a;b,c;d|q) = \phi_x(\ga;
\tilde a;\tilde b,\tilde c;\tilde d|q).
\tag\eq{140}
$$
Assume that (i) $0<b,c\leq a<d/q$, (ii)
$bd, cd\geq q$ and (iii) $ab,ac<1$, and let 
$t<0$ be an extra parameter.
By ${\Cal H}(a;b,c;d;t)$ we denote the weighted $L^2$-space
of symmetric functions $f(x)=f(x^{-1})$ 
with respect to the measure $d\nu(x;a,b,c;d|q,t)$, 
which is defined as follows. 
We introduce
$$
\De(x) = {{(x^{\pm 2}, qx^{\pm 1}/d;q)_\infty}\over{
\th(tdx^{\pm 1}) \, (ax^{\pm1}, bx^{\pm 1}, cx^{\pm 1};q)_\infty}},
\tag\eq{145}
$$
where $\th(x)=(x,q/x;q)_\infty$ is the 
(renormalised) Jacobi theta-function and where
we use $(cx^{\pm 1};q)_\infty=(cx,c/x;q)_\infty$ and
similarly for the other $\pm$-signs. 
Note that this the standard Askey-Wilson weight function, see
\cite{\AskeW}, \cite{\GaspR, Ch.~6}, multiplied by 
$q$-constant function that can be written as 
a quotient of theta-functions. The positive 
measure is now given by
$$
\int_{\C^\ast} f(x) d\nu(x;a;b,c;d|q,t) =
{K\over{4\pi i}}\int_{\T} f(x)\De(x){{dx}\over x} + 
K \sum_{s\in S} f(s) \text{Res}_{x=s}
{{\De(x)}\over x},
\tag\eq{150}
$$
where $S=S_+\cup S_-$, $S_+=\{ aq^k\mid k\in\Zp,\, aq^k>1\}$ and
$S_-=\{ tdq^k\mid k\in\Z,\, tdq^k<1\}$, under the
generic assumption on the parameters that $S\cup S^{-1}$
consists of simple poles of $\De$. Note that this
condition can be removed by extending the definition of the
masses at the discrete set by continuity, see \cite{\KoelSAW}
for details. 
Here $K=K(a;b,c;d;t)$ is a positive constant defined by
$$
K= (qabcdt^2)^{-\hf} (ab,ac,bc,qa/d,q;q)_\infty 
\bigl( \th(qt)\th(adt) \th(bdt) \th(cdt)\bigr)^\hf.
\tag\eq{155}
$$
The Askey-Wilson function transform pair for a
sufficiently nice 
function $u\in {\Cal H}(a;b,c;d;t)$ is given by 
$$
\aligned
\hat u(\ga) &= \int_{\C^\ast} u(x) \phi_\ga(x;a;b,c;d|q) 
\, d\nu(x;a;b,c;d|q,t), \\ 
u(x) &= \int_{\C^\ast} \hat u(\ga) \phi_\ga(x;a;b,c;d|q) 
\, d\nu(\ga;\tilde a;\tilde b,\tilde c; 
\tilde d|q, \tilde t),
\endaligned
\tag\eq{160}
$$
where $\tilde t=1/qadt$. Furthermore, $(\tilde a, \tilde b,\tilde
c,\tilde d,\tilde t)$ satisfy the same conditions as
$(a,b,c,d,t)$. 
So, in view of \eqtag{140}, 
the inverse of the Askey-Wilson function transform is
the Askey-Wilson transform for the dual set of parameters.
Moreover, the Askey-Wilson function transform 
extends to an isometric isomorphism from ${\Cal H}(a,b,c;d;t)$ 
to ${\Cal H}(\tilde a,\tilde b,\tilde c;\tilde d;\tilde t)$. 
Proofs of these results can be found in \cite{\KoelSAW}. See
Suslov \cite{\Susl}, \cite{\SuslPP} for Fourier-Bessel type
orthogonality relations for the Askey-Wilson functions. 

\subhead \the\sectionnumber.2. The Askey-Wilson 
polynomials\endsubhead The Askey-Wilson polynomials are 
the eigenfunctions of $L$ as in \eqtag{110}, \eqtag{120},
which are polynomial in $\hf(x+x^{-1})$. These 
orthogonal polynomials
are very well-known, see \cite{\AskeW}, \cite{\GaspR, \S 7.5},
and they are on top of the Askey scheme of basic hypergeometric
orthogonal polynomials, see \cite{\KoekS}. See also
Brown, Evans and Ismail \cite{\BrowEI} 
for an operator approach to the 
Askey-Wilson polynomials and their orthogonality relations. 
The Askey-Wilson functions reduce to the Askey-Wilson polynomials
for $\ga^{-1}=\tilde a q^n$, $n\in\Zp$, 
since the ${}_8W_7$-series in
\thetag{390} reduces to a terminating series 
$$
\multline
{}_8W_7(aq^{-n}/d;ax,a/x,q^{-n},q^{1-n}/cd,q^{1-n}/bd;q,q^nbc) = \\
{{(aq^{1-n}/d,q^{1-n}/ad;q)_n}\over{(q^{1-n}/dx, q^{1-n}x/d;q)_n}}
\, {}_4\vp_3\left( {{q^{-n},abcdq^{n-1},ax,a/x}\atop{
ab,ac,ad}};q,q\right)
\endmultline
\tag\eq{170}
$$
by \cite{\GaspR, (III.18)}. We can also use \cite{\GaspR, (III.36)}
to write the ${}_8W_7$-series as a sum of two balanced
${}_4\vp_3$-series, which reduces to a single terminating
balanced ${}_4\vp_3$-series for $\ga^{-1}=\tilde a q^n$, 
see e.g. Suslov \cite{\Susl}, \cite{\SuslPP}. This shows that the
Askey-Wilson function 
$\phi_\ga(x;a;b,c;d|q)$ is the analytic continuation of the
Askey-Wilson polynomial. 

\subhead \the\sectionnumber.3. The Askey-Wilson 
$q$-Bessel functions\endsubhead 
In \eqtag{130} we replace $c$ by $c\ep$, $d$ by $d/\ep$ and $\ga$ 
by $\ga\ep$, then the formal limit transition gives 
$$
\lim_{\ep\downarrow 0}
\phi_{\ga\ep}(x;a;b,c\ep;{d\over\ep}|q) = 
{}_2\vp_1\left( {{ax, a/x}\atop{ab}};q, {q\over{\tilde d \ga}}
\right).
\tag\eq{180}
$$
Taking the limit in the second order $q$-difference equation
shows that the Askey-Wilson $q$-Bessel function
$$
J_\ga(x;a;b|q) ={}_2\vp_1\left( {{ax,a/x}\atop{ab}};q, 
-{q\ga\over{a}}\right)
\tag\eq{185}
$$
is a solution to $LJ_\ga(\cdot;a;b|q)= \ga J_\ga(\cdot;a;b|q)$,
with $L$ as in \eqtag{110} with 
$$
A(x) = {{(1-ax)(1-bx)x}\over{(1-x^2)(1-qx^2)}}, \qquad 
B(x)=A(x^{-1}).
\tag\eq{187}
$$
Taking the limit transition \eqtag{180} 
through the sequence $\ep=q^m$, $m\to\infty$, 
in the Askey-Wilson function transform pair with $t$ replaced
$\ep t$ formally leads to the following
orthogonality relations
$$
\multline
\int_{\C^\ast} J_{\ga q^k}(x;a;b|q)\, J_{\ga q^l}
(x;a;b|q)\, d\nu(x;a;b,q\ga;\ga^{-1}|q,-1) = \\
\de_{k,l} a^{-2k} {{(-aq^{-k}/\ga;q)_\infty}\over
{(-bq^{-k}/\ga;q)_\infty}} {{K(a;b,q\ga;\ga^{-1};-1)}
\over{ \bigl( (ab;q)_\infty \th(-a/\ga)\bigr)^2}}, 
\qquad k,l\in\Z, 
\endmultline
\tag\eq{190}
$$
where $\ga=-a\tilde t$ 
for $a,b,\ga>0$, $ab<1$, $a>b$ and $K$ as in 
\eqtag{155}. Observe that cancellation in 
the weight $\De$ \eqtag{145}
occurs in \eqtag{190}, since $cd=q$. 
The orthogonality relations can be obtained
from Kakehi \cite{\Kake}, see also \cite{\KoelSsu, App.~A}
where also other ranges of the parameters are considered. 
Moreover, 
the functions $J_{\ga q^k}(\cdot;a;b|q)$ form a complete
set in the weighted $L^2$-space for the measure in 
\eqtag{190}. 
Note that under the limit transition \eqtag{180} only the
infinite set discrete mass points that tend to $-\infty$ in 
\eqtag{160} survive.  

Using the same limit transition \eqtag{180} in \eqtag{170} 
we find the Askey-Wilson $q$-Bessel functions
\eqtag{185} with the corresponding orthogonality 
relations \eqtag{190} from the orthogonality
relations of the
Askey-Wilson polynomials, see also \cite{\KoelPhD, Ch.~3}. 
The Askey-Wilson $q$-Bessel function is also studied by
Bustoz and Suslov \cite{\BustS}, who derive Fourier series
expansions, and by Ismail, Masson and Suslov \cite{\IsmaMS},
who derive Fourier-Bessel 
type orthogonality relations for $J_\ga(x;a;b|q)$. 

\head \newsec Big $q$-analogue of the Jacobi function 
scheme\endhead

In this section we consider the analogue of Figure~1.1 on the
level of the big $q$-Jacobi case. So we consider
the second order difference
operator $L$ as in \eqtag{110} with  
$$
A(x) = a^2(1+{1\over{abx}})(1+{1\over{acx}}),
\qquad B(x)= (1+{q\over{bcx}})(1+{1\over x}). 
\tag\eq{220}
$$
The general set of eigenfunctions for $L$ in \eqtag{110} with
$A$ and $B$ as in \eqtag{220} is given by 
Gupta, Ismail and Masson \cite{\GuptIM}. 

\subhead \the\sectionnumber.1. The big $q$-Jacobi 
functions\endsubhead 
The big $q$-Jacobi functions are defined by
$$
\phi_\ga(x;a;b,c;q) = {}_3\vp_2\left( {{a\ga,a/\ga,-1/x}\atop
{ab,ac}};q, -bcx\right), \qquad |bcx|<1,
\tag\eq{230}
$$
and they satisfy $L\phi_\ga(\cdot;a;b,c;q)=
(-1-a^2+a(\ga+\ga^{-1}))\phi_\ga(\cdot;a;b,c;q)$. 
We can extend the definition of the big $q$-Jacobi function
\eqtag{230} to generic values of $x$ by requiring that this
second order $q$-difference equation remains valid, see
\cite{\KoelSbig}. 
Assume that the parameters $a$, $b$ and $c$ are positive, 
$a$ greater than $b$ and $c$,  
and that all pairwise products are less than one, and let
$z>0$. 
The big $q$-Jacobi function transform is given by 
the following transform pair
$$
\aligned
\hat u(\ga) &= \int_{-1}^{\infty(z)} u(x) \phi_\ga(x;a;b,c;q) 
{{(-qx,-bcx;q)_\infty}\over{(-abx,-acx;q)_\infty}}\, d_qx ,\\
u(x) &= C \int_{\C^\ast} \hat u(\ga) \, \phi_\ga(x;a;b,c;q) \, 
(\ga^{\pm 1}abc;q)_\infty^{-1}\, 
d\nu(\ga;a;b,c;q/abc|q,-z^{-1}),\\
C &= {{\th(-abz)\th(-acz)\th(-bcz)\, (ab,ac;q)_\infty^2}
\over{(1-q)z\, \th(-qz) K(a;b,c;q/abc;-1/z)}}
\endaligned
\tag\eq{250}
$$
with the notation of \eqtag{150} and 
where the $q$-integral is defined by 
$$
\int_{-1}^{\infty(z)} f(x) \, d_qx = 
(1-q)\sum_{k=0}^\infty f(-q^k)q^k + (1-q)z\sum_{k=-\infty}^\infty
f(zq^k)q^k.
$$
After a suitable scaling the big $q$-Jacobi function
transform extends to an isometric isomorphism of 
the weighted $L^2$-space with respect to the $q$-integral
with weight as in the first equality of \eqtag{250} 
onto the weighted $L^2$-space of symmetric functions
with respect to the weight $C(\ga^{\pm 1}abc;q)_\infty^{-1}\, 
d\nu(\ga;a,b,c;q/abc|q,-1/z)$. 
See \cite{\KoelSbig} for a proof of these statements. 

\subhead \the\sectionnumber.2. The big $q$-Jacobi 
polynomials\endsubhead 
The polynomial eigenfunctions to \eqtag{110}, \eqtag{220}
are the big $q$-Jacobi polynomials and they 
are given by the big $q$-Jacobi functions \eqtag{230} with
$\ga=aq^n$, $n\in\Zp$, 
$$
{}_3\vp_2\left( {{a^2q^n,q^{-n},-1/x}
\atop{ab,ac}};q,-bcx\right) = 
{{(cq^{-n}/a;q)_n}\over{(ac;q)_n}}
{}_3\vp_2\left( {{q^{-n}, a^2q^n,-abx}\atop{ab,qa/c}};q,q\right),
\tag\eq{255}
$$
see \cite{\KoelSbig, Prop.~5.3}. The big $q$-Jacobi polynomials
are orthogonal with respect to a positive discrete measure
supported on $-q^{\Zp}\cup -q^{\N}/bc$ for 
$ab<1$, $ac<1$, $qa/b<1$, $qa/c<1$ and $bc<0$, see
\cite{\KoelSbig, \S 10}. See Andrews and Askey \cite{\AndrA},
or \cite{\GaspR, \S 7.3}, \cite{\KoekS} for the standard 
definition of the big $q$-Jacobi polynomials. 

\subhead \the\sectionnumber.3. The big $q$-Bessel  
functions\endsubhead If we replace $c$ by $c\ep$ and $\ga$ by
$\ga\ep$ in \eqtag{230}, we obtain the following limit
$$
\lim_{\ep\downarrow 0} \phi_{\ga\ep}(x;a;b,c\ep;q) 
= {}_1\vp_1\left( {{-1/x}\atop{ab}};q, -{{abcx}\over{\ga}}
\right).
\tag\eq{260}
$$
We define the big $q$-Bessel function by
$$
J_\ga(x;a;q) = {}_1\vp_1\left( {{-1/x}\atop{a}};q, a\ga x\right).
\tag\eq{265}
$$  
The big $q$-Bessel function
is a solution to $LJ_\ga(\cdot;a;q)=-\ga
J_\ga(\cdot;a;q)$, with $L$ as in \eqtag{110} with
$A(x)=x^{-1}(1+1/ax)$ and $B(x)=q(1+1/x)/ax$. Taking the
limit in the big $q$-Jacobi function transform pair
\eqtag{250} with $z$ fixed, shows that for $\ga>0$, $0<a<1$, 
$$
\int_{-1}^{\infty(q/a\ga)} 
\bigl( J_{\ga q^k}J_{\ga q^l}\bigr)(x;a;q)
{{(-qx;q)_\infty}\over{(-ax;q)_\infty}}\, d_qx = \de_{k,l} 
(1-q) {{(q;q)_\infty^2 \th(-a\ga)}\over{(a;q)_\infty^2 \th(-\ga)}}
a^{-k} (-q^k\ga;q)_\infty
\tag\eq{270}
$$
for $k,l\in\Z$. Here the extra parameter $z$ in the measure
of the big $q$-Jacobi function transform \eqtag{250} is 
in the limit inverse proportional to $\ga$. 
Under the limit transition \eqtag{260} the only part of the
spectrum of \eqtag{250} that survives, is the infinite set of
discrete mass points tending to $-\infty$. 
Moreover, the big $q$-Bessel functions 
$J_{\ga q^l}(\cdot;a;q)$ form a complete orthogonal set in the 
weighted $L^2$-space for the measure in \eqtag{270}. See 
\cite{\CiccKK} for the proof of \eqtag{270}.  

The big $q$-Bessel function in \eqtag{260} can also
be obtained by taking the limit from the big $q$-Jacobi polynomials
and then the orthogonality relations \eqtag{270} 
can be obtained from the orthogonality relations for the big
$q$-Jacobi polynomials in a rigorous way, see \cite{\CiccKK, \S 6}.

\head \newsec Little $q$-analogue of the Jacobi function 
scheme\endhead

In this section we consider the analogue of Figure~1.1 on the
level of the little $q$-Jacobi case. So the second order difference
operator $L$ is now as in \eqtag{110} with 
$$
A(x) = a^2(1+{1\over{ax}}),
\qquad B(x)= (1+{q\over{bx}}). 
\tag\eq{320}
$$
The general set of eigenfunctions for $L$ with
$A$ and $B$ as in \eqtag{320} is given by the solutions
to second order $q$-difference equation that is the $q$-analogue
of the hypergeometric differential equation, see \cite{\GaspR,
Ch.~1}.

\subhead \the\sectionnumber.1. The little $q$-Jacobi 
functions\endsubhead 
The little $q$-Jacobi functions are defined by
$$
\phi_\ga(x;a;b;q) = {}_2\vp_1\left( {{a\ga,a/\ga}\atop
{ab}};q, -bx\right),
\tag\eq{330}
$$
and they satisfy $L\phi_\ga(\cdot;a;b;q)=
(-1-a^2+a(\ga+\ga^{-1}))\phi_\ga(\cdot;a;b;q)$. For $y>0$,
$a>b>0$, and $ab<1$ we have the transform pair
$$
\aligned
\hat u(\ga) &= \sum_{k=-\infty}^\infty u(k) 
\phi_\ga(yq^k;a;b;q) \, a^{2k}{{(-q^{1-k}/ay;q)_\infty}
\over{(-q^{1-k}/by;q)_\infty}}, \\
u(k) &=  C \int_{\C^\ast} \hat u(\ga) \phi_\ga(yq^k;a;b;q) 
\, d\nu(\ga;a;b,aby;q/aby|q,-1),
\endaligned
\tag\eq{340}
$$
for $C=(ab;q)_\infty^2\th(-by)^2/K(a;b,aby;q/aby;-1)$ 
and using the
notation \eqtag{150}. 
Observe that cancellation in 
the weight $\De$ \eqtag{145}
occurs in \eqtag{340}, since $cd=q$. 
For a proof of \eqtag{340} see Kakehi \cite{\Kake} or 
\cite{\KoelSsu, App.~A}, where the result is given for 
more general parameter values. 

\subhead \the\sectionnumber.2. The little $q$-Jacobi 
polynomials\endsubhead The little $q$-Jacobi polynomials
are the polynomial eigenfunctions of $L$ as in \eqtag{110}, 
\eqtag{320}, and they occur for $\ga=aq^n$, $n\in\Zp$;
$$
\phi_{aq^n}(x;a;b;q)= {}_2\vp_1\left( {{q^{-n}, a^2q^n}\atop{ab}}
;q, -bx\right).
\tag\eq{350}
$$
This is not the standard expression for the little $q$-Jacobi
polynomials as introduced by Andrews and Askey \cite{\AndrA},
or see \cite{\GaspR, \S 7.3}, \cite{\KoekS}, where the
orthogonality relations can be found. 

\subhead \the\sectionnumber.3. The little $q$-Bessel 
functions\endsubhead If we replace $\ga$ by $\ga\ep$ and $x$ by
$x\ep$, we obtain the little $q$-Bessel functions from 
the little $q$-Jacobi functions;
$$
\lim_{\ep\downarrow 0} \phi_{\ga\ep}(x\ep;a;b;q) = 
{}_1\vp_1\left( {0\atop{ab}};q, -{{abx}\over{\ga}}\right). 
\tag\eq{360}
$$
The little $q$-Bessel function is defined by
$$
j_\ga(x;a;q)= {}_1\vp_1(0;a;q,q\ga x).
\tag\eq{365}
$$
Note that the little $q$-Bessel function is self-dual;
$j_\ga(x;a;q)=j_x(\ga;a;q)$. 
These $q$-analogues of the Bessel function are also known under
the name ${}_1\vp_1$ $q$-Bessel function or Hahn-Exton $q$-Bessel
function, see \cite{\KoorS} and \cite{\KoelJAT} for 
historic references to Hahn, Exton and Jackson. 
The little $q$-Bessel functions are eigenfunctions of
$Lj_\ga(\cdot;a;q)=-q\ga j_\ga(\cdot;a;q)$ with $L$ as in
\eqtag{110} with $A(x)=a/x$ and $B(x)=q/x$. This is the
second order $q$-difference equation for the ${}_1\vp_1$-series,
and we obtain this from the second order $q$-difference
equation for the little $q$-Jacobi function.  
The little $q$-Bessel functions satisfy the orthogonality
relations for $0<a<1$, see Koornwinder and 
Swarttouw \cite{\KoorS, Prop.~2.6},
$$
\sum_{k=-\infty}^\infty a^k j_{q^n}(q^k;a;q)\, 
j_{q^m}(q^k;a;q) = \de_{n,m} a^{-n}
{{(q;q)_\infty^2}\over{(a;q)_\infty^2}}, \qquad n,m\in\Z.
\tag\eq{370}
$$
In the limit transition \eqtag{360} of the little $q$-Jacobi
function transform \eqtag{340} the only part of the spectrum
that survives the contraction is the infinite set of discrete
mass points tending to $-\infty$, which leads to a formal
derivation of \eqtag{370}. We note that the extra degree of
freedom in the measure of the little $q$-Jacobi function
transform \eqtag{340} drops out in the limit. 
By self-duality we see that the little $q$-Bessel functions
form a complete orthogonal set with respect to the discrete
measure in \eqtag{370}. 

The same limit transition of the little $q$-Jacobi polynomials
to the little $q$-Bessel functions as in \eqtag{360} is 
valid, and in the limit the orthogonality relations for the
little $q$-Jacobi polynomials tend to the orthogonality
relations \eqtag{370}, see \cite{\KoorS} for a rigorous
proof. 

\head \newsec Limit transitions\endhead

In the previous sections limits within each level of 
the $q$-analogue of the Jacobi function schemes have been
discussed. Now there are also limits from the Askey-Wilson 
polynomials to the big $q$-Jacobi polynomials and from the big
$q$-Jacobi polynomials to the little $q$-Jacobi polynomials,
see \cite{\KoorSIAM}, 
\cite{\StokK} or \cite{\KoekS}. 
In this section we show that these limit transitions also
hold for the appropriate analogues of the Jacobi and Bessel
function. 

\subhead \the\sectionnumber.1. Limit from the Askey-Wilson
case to the big $q$-case\endsubhead 
In the Askey-Wilson function of \eqtag{130} we replace
$(a,b,c,d,x)$ by $(a/\ep,b\ep,c\ep,d/\ep,-x/\ep)$. Then,
using \eqtag{132}, the Askey-Wilson function tends to the
big $q$-Jacobi function \eqtag{230} as $\ep\downarrow 0$; 
$$
\multline
\lim_{\ep\downarrow 0} 
\phi_\ga(-{x\over\ep};\frac{a}{\ep};b\ep,c\ep,{d\over\ep}|q) 
= {{(-q\tilde ax\ga/d,\tilde c/\ga;q)_\infty}
\over{(\ga\tilde c,ac, qa/d, -qx/d;q)_\infty}} 
\, {}_3\vp_2\left( {{-bx,\tilde a\ga,\tilde b\ga}\atop
{ab,-q\tilde a\ga x/d}};q,{{\tilde c}\over{\ga}}\right) \\
= \frac{1}{(qa/d;q)_\infty} \phi_\ga(\frac{x}{a};\tilde a;
\tilde b,\tilde c;q),
\endmultline
\tag\eq{410}
$$
where we have used \cite{\GaspR, (III.9)} in the second equality.
Keeping $t$ fixed we can formally take the limit
transition within the Askey-Wilson function transform 
pair \eqtag{160} to recover the big $q$-Jacobi function transform
\eqtag{250} with parameters $(a,b,c,z)$ replaced by 
$(\tilde a, \tilde b, \tilde c, -td/a)$. 
Taking $\ga=\tilde a q^n$, $n\in\Zp$,
gives back the limit transition from the Askey-Wilson
polynomials to the big $q$-Jacobi polynomials. 

In the Askey-Wilson $q$-Bessel function we replace
$(a,b,x,\ga)$ by $(a/\ep, b\ep,-x/\ep,\ga\ep)$ in \eqtag{185}
and take the limit $\ep\downarrow 0$, which gives the 
big $q$-Bessel function;
$$
\lim_{\ep\downarrow 0} 
J_{\ga\ep}(-{x\over\ep};\frac{a}{\ep};b\ep|q)= 
J_{q\ga/b}(\frac{x}{a};ab;q).
\tag\eq{420}
$$
In this limit transition the second order $q$-difference equation
for the Askey-Wilson $q$-Bessel function goes over into
the second order $q$-difference equation
for the big $q$-Bessel function. 
And the orthogonality relations \eqtag{190} go over into the
orthogonality relations \eqtag{270}, since only the 
discrete mass points of the measure survive in the limit.

\subhead \the\sectionnumber.2. Limit from the big $q$-case to 
the little $q$-case\endsubhead 
The big $q$-Jacobi function of \eqtag{230} tends to the 
little $q$-Jacobi function in \eqtag{330} by
$$
\lim_{c\downarrow 0} 
\phi_\ga({x\over c};a;b,c;q) = \phi_\ga(x;a;b;q).
\tag\eq{430}
$$
In this limit transition the big $q$-Jacobi function transform
\eqtag{250} tends to the little $q$-Jacobi function after
taking $z=y/c$ for the extra parameter in the big $q$-Jacobi
function transform \eqtag{250}. The
second order $q$-difference equation for the big $q$-Jacobi
functions tends to the second order $q$-difference equation for
the little $q$-Jacobi functions under \eqtag{430}. 

In the big $q$-Bessel function \eqtag{265} we replace
$x$ by $x/c$ and $\ga$ by $cq\ga/a$ and take the limit
$$
\lim_{c\downarrow 0} J_{cq\ga/a}({x\over c};a;q) = j_\ga(x;a;q)
\tag\eq{440}
$$
to find the little $q$-Bessel function \eqtag{365}.
In this limit transition the big $q$-Hankel orthogonality
relations \eqtag{270} tend to the little $q$-Hankel  
orthogonality relations \eqtag{370}, for which we observe that
the $\ga$-dependence drops out in the limit.  
The second order $q$-difference equation for the 
big $q$-Bessel functions tends to the second order 
$q$-difference equation for the little $q$-Bessel functions
under \eqtag{440}.
 
\head \newsec Duality and factoring of limits\endhead

As already observed, the Askey-Wilson function transform and
the little $q$-Hankel transform are self-dual. However, the
transforms in between are not self-dual, and the dual transforms
are described in this section. It turns out that these are 
related to indeterminate moment problems, see \cite{\CiccKK}, 
\cite{\KoelSbig}. 
Moreover, we find that the
limit transition of the little $q$-Jacobi function to the little
$q$-Hankel transform and of the Askey-Wilson function to the
Askey-Wilson $q$-Hankel transform factors through the transforms
dual to the big $q$-Hankel transform and dual 
to the big $q$-Jacobi function transform. 

\subhead \the\sectionnumber.1. Duality between the little $q$-Jacobi
function transform and the Askey-Wilson $q$-Hankel 
transform\endsubhead 
The little $q$-Jacobi function of \eqtag{330} and the
Askey-Wilson $q$-Bessel function of \eqtag{185} are related
by interchanging the geometric parameter $x$ and the spectral
parameter $\ga$;
$$
J_\ga(x;a;b|q) = \phi_x({{q\ga}\over{ab}};a,b;q).
\tag\eq{510}
$$
The orthogonality relations \eqtag{190} 
and the transform \eqtag{340} can be obtained from each other
by using \eqtag{510}. So the little $q$-Jacobi
function and the Askey-Wilson $q$-Bessel function 
are also eigenfunctions of  
a second order $q$-difference equation acting on the 
spectral parameter. 

\subhead \the\sectionnumber.2. The dual to the big $q$-Hankel
transform\endsubhead 
The big $q$-Bessel functions form an orthogonal basis for 
the weighted $L^2$-space on $[-1,\infty(q/a\ga))_q$ 
described in \eqtag{270},  see
\cite{\CiccKK}. The corresponding dual orthogonality relations
are then labeled by the support of this measure, i.e. 
by $-q^{\Zp}$ and by $q^\Z/a\ga$. In the first case, the 
big $q$-Bessel functions at $x=-q^n$ are related to 
the $q$-Laguerre polynomials in $\ga$;
$$
J_\ga(-q^n;a;q)= {}_1\vp_1\left( {{q^{-n}}\atop{a}};q,-a\ga q^n
\right),
\tag\eq{520}
$$
and the dual orthogonality relations for $J_\ga(-q^n;a;q)$ and
$J_\ga(-q^m;a;q)$, $n,m\in\Zp$, reduce
to the orthogonality relations for the $q$-Laguerre
polynomials related to Ramanujan's ${}_1\psi_1$-sum, see
Moak \cite{\Moak, Thm.~2}, or \cite{\GaspR, Exer.~7.43(ii)}. 
The big $q$-Bessel functions 
are eigenfunctions of 
a second order $q$-difference operator in the spectral parameter
$\ga$; 
$$
\Bigl( (1+{1\over\ga}) (T_q^\ga -1) +
{q\over{a\ga}}(T_{q^{-1}}^\ga-1)\Bigr) J_\ga(x;a;q) = 
-(1+x) J_\ga(x;a;q), 
\tag\eq{540}
$$
where $(T_{q^{\pm 1}}^\ga f)(\ga)= f(q^{\pm 1}\ga)$. Note that
\eqtag{540} is nothing but the second order $q$-difference
equation for the ${}_1\vp_1$-series. 

It is known that the $q$-Laguerre polynomials correspond
to an indeterminate moment problem and that this solution
to the moment problem is not extremal in the sense of
Nevannlina, meaning that the $q$-Laguerre polynomials are
not dense in the corresponding weighted $L^2$-space. 
The functions of $\ga$ defined by, $p\in\Z$,  
$$
J_\ga({{q^{p+1}}\over{a\ga}};a;q) = 
{}_1\vp_1\left( {{-a\ga q^{-1-p}}\atop a};q, q^{p+1}\right) 
= {{(q^{p+1};q)_\infty}\over{(a;q)_\infty}}
\, {}_1\vp_1\left( {{-\ga}\atop{ q^{p+1}}};q,a\right)
\tag\eq{530}
$$
display $q$-Bessel coefficient behaviour. Here we used a
limit case of the transformation formula \cite{\GaspR, (III.2)}. 
These $q$-Bessel coefficients complement the orthogonal 
set of $q$-Laguerre polynomials into an orthogonal basis of
the corresponding weighted $L^2$-space. This is
a direct consequence of of the orthogonality
relations \eqtag{270} and the completeness. 
See  \cite{\CiccKK} for
details, where also the spectral analysis of \eqtag{540}
is given. 

\subhead \the\sectionnumber.3. The dual to the big $q$-Jacobi
function transform\endsubhead Evaluating the big $q$-Jacobi
function at the point $-q^k$, $k\in\Zp$, gives a terminating
series in which the base is inverted to $q^{-1}$;
$$
\phi_\ga(-q^k;a;b,c;q) = {}_3\vp_2\left( {{q^k, \ga/a,1/a\ga}
\atop{1/ab,1/ac}};q^{-1},q^{-1}\right).
\tag\eq{570}
$$
The right hand side is a polynomial of degree $k$ in
$\hf(\ga+\ga^{-1})$, which is a continuous dual $q^{-1}$-Hahn 
polynomial with parameters $a^{-1}$, $b^{-1}$, $c^{-1}$, 
i.e. an Askey-Wilson polynomial with one parameter set to zero,  
see \cite{\GaspR, \S 7.5}, \cite{\KoekS}. Hence, the transform
\eqtag{250} gives us an orthogonality measure for the
continuous dual $q^{-1}$-Hahn polynomials, together with
a complementing set of orthogonal functions in 
$\hf(\ga+\ga^{-1})$, namely $\phi_\ga(zq^l;a;b,c;q)$, 
$l\in\Z$. See \cite{\KoelSbig, \S 9} for more details.
In particular we find that the big $q$-Jacobi functions are
eigenfunctions  to 
a second order $q$-difference operator in the spectral
parameter $\ga$, see \cite{\IsmaR}; 
$$
\gathered
\bigl( A(\ga)(T_q^\ga-1)+A(\ga^{-1})(T_{q^{-1}}^\ga-1)\bigr)
\phi_\ga(x;a;b,c;q)= -(1+x) \phi_\ga(x;a;b,c;q), \\
A(\ga)= {{(1-1/\ga a)(1-1/\ga b)(1-1/\ga c)}\over
{(1-\ga^{-2})(1-1/\ga^2 q)}},
\endgathered
\tag\eq{580}
$$
cf. the second order $q$-difference equation for the
continuous dual $q^{-1}$-Hahn polynomials \cite{\AskeW},
\cite{\KoekS}. 
See also Rosengren \cite{\Rose, \S 4.6} for another orthogonality
measure for the continuous dual $q^{-1}$-Hahn polynomials.

\subhead \the\sectionnumber.4. Remaining limit 
transitions\endsubhead 
Using duality and the limit transitions in \S 5 we obtain the
remaining limit transitions in Figure~1.2, namely the limits from
the Askey-Wilson functions to the continuous dual $q^{-1}$-Hahn
polynomials and the associated $q$-Bessel functions, and from
this family to the Askey-Wilson $q$-Bessel functions, and from the
little $q$-Jacobi functions to the $q$-Laguerre polynomials
and the associated $q$-Bessel functions, and from this family 
to the little $q$-Bessel functions. These limits are formal
and can be taken in the second order $q$-difference equation and
in the transform pairs. 

As an example, we illustrate the limit from the Askey-Wilson 
functions to the continuous dual $q^{-1}$-Hahn polynomials
and the associated $q$-Bessel functions. Using the duality 
\eqtag{140} we have for $x_k=-aq^k$, $k\in\Zp$, and 
$\ep>0$ sufficiently small,
$$
\phi_{-x_k/\ep}(\ga;\tilde a;\tilde b,\tilde c;
\frac{\tilde d}{\ep^2}|q)
= \phi_\ga(-\frac{x_k}{\ep};\frac{a}{\ep};b\ep,
c\ep;\frac{d}{\ep}|q).
\tag\eq{581}
$$
By \eqtag{410} the limit $\ep\downarrow 0$ gives 
$$
\lim_{\ep\downarrow 0} 
\phi_{-x_k/\ep}(\ga;\tilde a;\tilde b,\tilde c;
\frac{\tilde d}{\ep^2}|q) = \frac{1}{(qa/d;q)_\infty}
\phi_\ga(-q^k;\tilde a;\tilde b,\tilde c;q)
\tag\eq{582}
$$
and the right hand side is a continuous dual $q^{-1}$-Hahn
polynomial of degree $k$ up to a constant, see \eqtag{570}.
Replacing $x_k$ by the discrete weights $y_l=-tdq^l$, $l\in\Z$,
in \eqtag{581}, \eqtag{582} we see that 
$\phi_{-y_l/\ep}(\ga;\tilde a;\tilde b,\tilde c;
\frac{\tilde d}{\ep^2}|q)$ tend to $q$-Bessel coefficient type
functions which complement the continuous dual $q^{-1}$-Hahn
polynomials to an orthogonal basis in the corresponding
weighted $L^2$-space.  The non-extremal measure, which is
parametrised by the $z$-parameter in the big $q$-Jacobi function
transform \eqtag{250}, corresponds to $z=-td/a$. In this 
limit the Askey-Wilson function transform formally tends to the
the orthogonality relations for the continuous dual $q^{-1}$-Hahn
polynomials and the corresponding $q$-Bessel functions. 
For this we only need to remark that, by duality, this reduces to
limit transition \eqtag{410} 
of the Askey-Wilson function transform 
to the big $q$-Jacobi function transform. 

All the other cases can be considered in a similar manner
and for completeness we give the underlying limit transitions;
$$
\align
&\lim_{c\downarrow 0}\phi_{\ga}(\frac{x}{c};a;b,c;q)=
J_{abx/q}(\ga;a;b;q),
\tag\eq{5100}\\
&\lim_{\ep\downarrow 0} \phi_{\ga/\ep}(\ep x; \frac{a}{\ep}; 
\ep b;q) 
= J_{ax}(-\frac{\ga}{a};ab;q),
\tag\eq{550}\\
&\lim_{\ep\downarrow 0} J_{\ep\ga}({x\over\ep};a;q)= 
j_x({{\ga a}\over q};a;q).
\tag\eq{560}
\endalign
$$
Let us finally note that this makes the scheme of limit transitions
in Figure~1.2 into a commutative diagram. 

\head \newsec Concluding remarks\endhead

\subhead \the\sectionnumber.1. Quantum group theoretic
interpretation\endsubhead 
There is a group theoretic interpretation for the scheme of
Figure~1.1, see Koornwinder \cite{\KoorLNM}. Here the Jacobi
polynomials have an interpretation on the compact real Lie group 
$SU(2)$ as matrix elements of finite-dimensional 
irreducible unitary representations,
whereas the Jacobi functions have an interpretation as
matrix elements of infinite-dimensional 
irreducible unitary representations of the non-compact real 
Lie group $SU(1,1)$. The real Lie groups $SU(2)$ and $SU(1,1)$
are both real forms of the complex Lie group $SL(2,\C)$. 
Bessel functions occur as matrix elements of infinite-dimensional 
irreducible unitary representations of the group $E(2)$ of
motions of the Euclidean plane, and then the limit transition
can be interpreted on the level of Lie groups as a contraction.

There is also a quantum group theoretic interpretation for 
the Askey-Wilson function scheme in Figure~1.2. 
Here the Askey-Wilson, big and
little $q$-Jacobi polynomials are interpreted as matrix elements
of irreducible unitary representations of the quantum $SU(2)$ 
group, see e.g. \cite{\KoelAAM}, \cite{\KoorSIAM}, \cite{\Noum}, 
\cite{\NoumM} and further references. These interpretations
also naturally lead to the limit transitions from Askey-Wilson
polynomials to big $q$-Jacobi polynomials, and from 
big $q$-Jacobi polynomials to little $q$-Jacobi polynomials,
see Koornwinder \cite{\KoorSIAM}. 
The various $q$-Hankel transforms
and $q$-Bessel functions 
have an interpretation on the quantum $E(2)$ group (or better, its
twofold covering); see \cite{\KoelDMJ}, \cite{\VaksKmotion} for
the little $q$-Bessel function; see 
\cite{\BoneCGST}, \cite{\KoelJCAM} 
for the big $q$-Bessel function, and see \cite{\KoelPhD},
\cite{\KoelIM} for the Askey-Wilson $q$-Bessel function. The 
limit transitions from the $q$-analogues of the Jacobi polynomials 
to the $q$-analogues of the Bessel functions as described in
\S\S 2,3, 4 
are motivated from a similar
contraction from the quantum $SU(2)$ group to the quantum $E(2)$
group. The quantum $SU(2)$ group is a real form of the quantum
$SL(2,\C)$ group, and one of the other real forms is the
quantum $SU(1,1)$ group. The interpretation 
of the little $q$-Jacobi
functions on the quantum $SU(1,1)$ group as matrix elements
of irreducible unitary representations 
is due to 
Masuda et al. \cite{\MasuMNNSU}, Kakehi \cite{\Kake}, 
Kakehi, Masuda and Ueno \cite{\KakeMU} and Vaksman and 
Korogodski\u\i\ \cite{\VaksKsu}. For the big $q$-Jacobi 
and Askey-Wilson function such an interpretation 
is given in \cite{\KoelSsu}. This interpretation 
leads in a natural way to the limit transitions
of the Askey-Wilson functions to big $q$-Jacobi functions, and
from big $q$-Jacobi functions to little $q$-Jacobi functions
as described in \S 5. 
Moreover, from a contraction procedure 
from the quantum $SU(1,1)$ group to the quantum $E(2)$ group
we obtain the limit transition of the $q$-analogues of the
Jacobi function to the corresponding $q$-analogues of the Bessel
function, see \S\S 2, 3 and 4. 
Due to these representation theoretic interpretation of the
three types of $q$-special functions on the quantum $SU(1,1)$,
$E(2)$ and $SU(2)$ groups, we may view them
as $q$-analogues of the Jacobi functions, Bessel functions and
Jacobi polynomials, respectively. This explains the naming of
the $q$-special functions in Figure~1.2.

The second order difference equation for the $q$-special
functions follows from the action of the Casimir element
for these quantum groups. In  \S 6.1 we have seen 
that the Askey-Wilson $q$-Bessel functions and the little
$q$-Jacobi functions both satisfy a second order $q$-difference
equation in the spectral parameter. The spectral parameter
corresponds to the representation label of the irreducible
representations of the quantum $E(2)$ group and 
the quantum $SU(1,1)$ group respectively. 
The second order $q$-difference equation in the 
spectral parameter can then be obtained from the tensor
product decomposition of a three-dimensional (non-unitary)
representation with an irreducible infinite 
dimensional representation into three irreducible infinite
dimensional representations of the corrresponding quantum groups,
see e.g. \cite{\KoelAAM, Remark~7.2} for the corresponding
statement for the quantum $SU(2)$ group. 

\subhead \the\sectionnumber.2. Further extensions\endsubhead
We briefly sketch some possible directions related to the
Askey-Wilson function scheme. 

{\bf 1.} From the previous subsection we see that the scheme 
depicted in Figure~1.2 is motivated by the simplest quantum groups;
namely for $SU(2)$, $SU(1,1)$, and $E(2)$. We would expect that
a greater extension of the scheme in Figure~1.2 is possible
using more complicated quantum groups, especially higher rank 
quantum groups, or other interpretations of special functions
on quantum groups, such as Clebsch-Gordan, Racah 
or other type of overlap coefficients.
See e.g. Koornwinder \cite{\KoorLNM} 
for a further extension of the
Jacobi function scheme in Figure~1.1.
In such an extension of Figure~1.2 the other $q$-analogues of
the Bessel function, as studied by Ismail \cite{\Isma}, 
should also find a place.

{\bf 2.} For the Hankel transform and Jacobi function transform
there are systems of orthogonal polynomials that 
are mapped onto each other, such as the Laguerre polynomials
for the Hankel transform or the Jacobi polynomials that are mapped
onto the Wilson polynomials by the Jacobi function transform,
see \cite{\KoorLNM}. It would be interesting to know what the 
corresponding results for the $q$-analogues of these 
transforms are.



\Refs

\ref\no \AndrA
\by G.E.~Andrews and R.~Askey
\paper Classical orthogonal polynomials
\inbook ``Polyn\^omes Orthogonaux et Applications''
\bookinfo LNM \vol 1171
\publaddr  Springer Verlag
\yr 1985
\pages 36--62
\endref

\ref\no \AskeW
\by R.~Askey and J.~Wilson
\paper Some basic hypergeometric orthogonal polynomials that
generalize Jacobi polynomials
\jour Mem. Amer. Math. Soc.
\vol 54
\issue 319
\yr 1985
\endref

\ref\no \BoneCGST
\by F.~Bonechi, N.~Ciccoli, R.~Giachetti, E.~Sorace, 
and M.~Tarlini
\paper Free $q$-Schr\"odinger equation from
homogeneous spaces of the $2$-dim Euclidean quantum group
\jour Comm. Math. Phys. \vol 175 \yr 1996
\pages 161--176
\endref

\ref\no \BrowEI
\by B.M.~Brown, W.D.~Evans and M.E.H.~Ismail 
\paper The Askey-Wilson polynomials and $q$-Sturm-Liouville 
problems 
\jour Math. Proc. Cambridge Philos. Soc. \bf 119 \yr 1996
\pages  1--16
\endref

\ref\no \BustS
\by J.~Bustoz and S.K.~Suslov 
\paper Basic analog of Fourier series on a $q$-quadratic grid
\jour Methods Appl. Anal. \vol 5 \yr 1998
\pages 1--38
\endref

\ref\no \CiccKK
\by N.~Ciccoli, E.~Koelink and T.H.~Koornwinder
\paper $q$-Laguerre polynomials and big $q$-Bessel functions
and their orthogonality relations
\jour Methods Appl. Anal. \vol 6 \yr 1999
\pages 109--127 
\endref

\ref\no \GaspR
\by G.~Gasper and M.~Rahman
\book Basic Hypergeometric Series
\publaddr Cambridge Univ. Press
\yr 1990
\endref

\ref\no \GuptIM 
\by D.P.~Gupta, M.E.H.~Ismail and D.R.~Masson
\paper Contiguous relations, basic hypergeometric functions, and
orthogonal polynomials III. Associated continuous 
dual $q$-Hahn polynomials
\jour J. Comput. Appl. Math. \vol 68 \yr 1996
\pages 115--149 
\endref

\ref\no \Isma
\by M.E.H.~Ismail
\paper The zeros of basic Bessel functions, the function
$J_{\nu+ax}(x)$, and associated orthogonal polynomials
\jour J. Math. Anal. Appl. \vol 86 \yr 1982
\pages 1--19
\endref

\ref\no \IsmaMS
\by M.E.H.~Ismail, D.R.~Masson and S.K.~Suslov
\paper The $q$-Bessel function on a $q$-quadratic grid
\inbook ``Algebraic Methods and $q$-Special Functions''
\bookinfo CRM Proc. Lect. Notes \vol 22
\eds J.F. van Diejen and L. Vinet
\publaddr AMS
\yr 1999
\pages 183--200
\endref

\ref\no \IsmaR 
\by M.E.H.~Ismail and M.~Rahman
\paper The associated Askey-Wilson polynomials 
\jour Trans. Amer. Math. Soc. \vol 328 \yr 1991
\pages 201--237
\endref
 
\ref\no \Kake
\by T.~Kakehi
\paper Eigenfunction expansion associated with the Casimir
operator on the quantum group $SU_q(1,1)$
\jour Duke Math. J.
\vol 80
\yr 1995 
\pages 535--573
\endref

\ref\no \KakeMU
\by T.~Kakehi, T.~Masuda and K.~Ueno
\paper Spectral analysis of a $q$-difference operator which 
arises from the quantum $SU(1,1)$ group
\jour J. Operator Theory
\vol 33
\yr 1995
\pages 159--196
\endref

\ref\no \KoekS
\by R.~Koekoek and R.F.~Swarttouw
\paper The Askey-scheme of hypergeometric orthogonal polynomials 
and its $q$-analogue
\paperinfo Delft University of Technology Report no. 98-17
\yr 1998
\endref

\ref\no \KoelPhD
\by H.T.~Koelink
\book On quantum groups and $q$-special functions
\bookinfo dissertation, Rijksuniversiteit Leiden
\yr 1991
\endref

\ref\no \KoelDMJ
\by H.T.~Koelink
\paper The quantum group of plane motions and the 
Hahn-Exton $q$-Bessel function \jour Duke Math. J. \vol 76 
\yr 1994 \pages 483--508
\endref

\ref\no \KoelIM
\by H.T.~Koelink
\paper The quantum group of plane motions and basic 
Bessel functions \jour Indag. Math. (N.S.)
\vol 6 \yr 1995 \pages 197--211
\endref

\ref\no \KoelJCAM
\by H.T.~Koelink
\paper Yet another basic analogue of Graf's addition formula
\jour J. Comput. Appl. Math. \vol 68
\yr 1996 \pages 209--220
\endref

\ref\no \KoelAAM
\by H.T.~Koelink
\paper Askey-Wilson polynomials and the
quantum $SU(2)$ group: survey and applications
\jour Acta Appl. Math.
\vol 44
\yr 1996
\pages 295--352
\endref

\ref\no \KoelJAT
\by H.T.~Koelink
\paper Some basic Lommel polynomials
\jour J. Approx. Theory \vol 96 \yr 1999
\pages 345--365
\endref

\ref\no \KoelSbig
\by E.~Koelink and J.V.~Stokman
\paper The big $q$-Jacobi function transform
\paperinfo preprint 40~p., {\tt math.CA/9904111}
\yr 1999
\endref

\ref\no \KoelSsu
\by E.~Koelink and J.V.~Stokman, with an appendix by M.~Rahman
\paper Fourier transforms on the quantum $SU(1,1)$ group 
\paperinfo preprint 77~p., {\tt math.QA/9911163} 
\yr 1999
\endref

\ref\no \KoelSAW
\by E.~Koelink and J.V.~Stokman
\paper The Askey-Wilson function transform
\paperinfo preprint
\toappear
\endref

\ref\no \KoorJF
\by T.H.~Koornwinder
\paper Jacobi functions and analysis on noncompact 
semisimple Lie groups
\inbook ``Special Functions: Group
Theoretical Aspects and Applications''
\eds R. A. Askey, T. H. Koornwinder and W. Schempp
\pages 1--85
\publaddr Reidel
\yr 1984
\endref

\ref\no \KoorLNM
\by T.H.~Koornwinder
\paper Group theoretic interpretations of Askey's scheme 
of hypergeometric orthogonal polynomials 
\inbook ``Orthogonal Polynomials and their Applications''
\pages 46--72 \bookinfo LNM \vol 1329
\publaddr Springer \yr 1988
\endref

\ref\no \KoorSIAM
\by T.H.~Koornwinder
\paper Askey-Wilson polynomials
as zonal spherical functions on the $SU(2)$ quantum group
\jour SIAM J. Math. Anal.
\vol 24
\yr 1993
\pages 795--813
\endref 

\ref\no \KoorLNSFQG 
\by T.H.~Koornwinder 
\paper Compact quantum groups and $q$-special functions
\inbook ``Representations of Lie
Groups and Quantum Groups'' \pages 46--128
\bookinfo Pitman Res. Notes Math. Ser. \vol 311
\publaddr Longman Sci. Tech. \yr 1994
\endref 

\ref\no \KoorS
\by T.H.~Koornwinder and R.F.~Swarttouw
\paper On $q$-analogues of the Fourier and Hankel transforms
\jour Trans. Amer. Math. Soc. \vol 333 \yr 1992 \pages 445--461
\endref

\ref\no \MasuMNNSU
\by T.~Masuda, K.~Mimachi, Y.~Nakagami, M.~Noumi, Y.~Saburi
and K.~Ueno
\paper Unitary representations of the quantum group $SU_q(1,1)$:
Structure of the dual space of $U_q(sl(2))$
\jour Lett. Math. Phys.
\vol 19
\yr 1990
\pages 197--194
\moreref
\paper II: Matrix elements of unitary representations and the
basic hypergeometric functions
\pages 195--204
\endref

\ref\no\Moak
\by D.S.~Moak
\paper The $q$-analogue of the Laguerre polynomials
\jour J. Math. Anal. Appl. \vol 81 \yr 1981
\pages 20--47
\endref

\ref\no \Noum
\by M.~Noumi
\paper Quantum groups and $q$-orthogonal polynomials. Towards a
realization of Askey-Wilson polynomials on $SU_q(2)$
\inbook ``Special Functions''
\eds M.~Kashiwara and T.~Miwa
\bookinfo ICM-90 Satellite Conference Proceedings
\publaddr Springer-Verlag
\yr 1991
\pages 260--288
\endref

\ref\no \NoumM
\by M.~Noumi and K.~Mimachi
\paper Askey-Wilson polynomials as spherical functions on $SU_q(2)$
\inbook ``Quantum Groups''
\bookinfo LNM \vol 1510
\publaddr Springer Verlag
\yr 1992
\pages 98--103
\endref

\ref\no \Rose
\by H.~Rosengren
\paper A new quantum algebraic interpretation of the
Askey-Wilson polynomials
\jour Contemp. Math.
\toappear
\endref

\ref\no \StokK 
\by J.V.~Stokman and T.H.~Koornwinder
\paper On some limit cases of Askey-Wilson polynomials
\jour J. Approx. Theory \vol 95 \yr 1998
\pages 310--330
\endref

\ref\no \Susl
\by S.K.~Suslov
\paper Some orthogonal very well poised ${}_8\phi_7$-functions
\jour J. Phys. A \vol 30 \yr 1997
\pages 5877--5885
\endref

\ref\no \SuslPP
\by S.K.~Suslov
\paper Some orthogonal very-well-poised 
${}_8\vp_7$-functions that generalize Askey-Wilson polynomials
\paperinfo preprint
\yr 1997
\endref

\ref\no \VaksKmotion
\by L.L.~Vaksman and L.I.~Korogodski\u\i
\paper The algebra of bounded functions on the quantum group 
of motions of the plane and $q$-analogues of Bessel functions
\jour Soviet Math. Dokl. \vol 39 \yr 1989 \pages 173--177
\endref

\ref\no \VaksKsu
\by L.L.~Vaksman and L.I.~Korogodski\u\i
\paper Spherical functions on the quantum group 
$\text{SU}(1,1)$ and a $q$-analogue of the
Mehler-Fock formula
\jour Funct. Anal. Appl. \vol 25 \yr 1991
\pages 48--49 
\endref

\ref\no \Wats
\by G.N.~Watson
\book A Treatise on the Theory of Bessel Functions
\publaddr Cambridge Univ. Press \yr 1944
\endref

\endRefs
\enddocument